# An arguable inconsistency in ZF

Bhupinder Singh Anand[1]

Classical theory proves that every primitive recursive function is strongly representable in PA; that PA and PRA can both be interpreted in ZF; and that if ZF is consistent, then PA+PRA is consistent. We show that PA+PRA is inconsistent; it follows that ZF, too, is inconsistent.

## 1. Overview

Classical theory proves that:

*a*) formal Peano Arithmetic, PA, and formal primitive recursive arithmetic, PRA, can both be interpreted in Zermelo-Fraenkel Set Theory, ZF.

*b*) if ZF is consistent, then PA+PRA is consistent.

*c*) every primitive recursive function is strongly representable in PA.

It also seems reasonable to suspect that:

*d*) every primitive recursive function cannot be defined in PA.

In the appended Meta-theorem 1 and Meta-lemma 1 [cf. An02], we now argue that:

*e*) PA+PRA is inconsistent.

---

[1] The author is an independent scholar. E-mail: re@alixcomsi.com; anandb@vsnl.com. Postal address: 32, Agarwal House, D Road, Churchgate, Mumbai - 400 020, INDIA. Tel: +91 (22) 2281 3353. Fax: +91 (22) 2209 5091.

In Meta-theorem 1, we consider Gödel's primitive recursive relation $\sim xB_{PA}(Sb(y\ 19|Z(y)))$ (cf. [Go31a], p24, def. 8.1; we note, in particular, that the index serves as a reminder that the primitive recursive relation $B_{PA}$ is specific to, and defined with reference to, the axioms of the recursively enumerable formal system PA, which is, essentially, Gödel's arithmetic P in [Go31a]), and argue that it is not the standard interpretation of any of its formal representations in PA.

In other words we argue that, if we assume the primitive recursive relation $\sim xB_{PA}(Sb(y\ 19|Z(y)))$ to be an abbreviation of some formula of PA, then we arrive at an inconsistency in PA.

Intuitively, this is unexceptionable, since, prima facie, the use of a representation theorem (cf. [Me64], Proposition 3.23, p131) in standard derivations (cf. [Me64], Proposition 3.31, p143) of Gödel's Theorem VI (cf. [Go31a], p24) appears to be both critical and necessary. Meta-theorem 1 proves the necessity, confirming (*d*).

We consequently note, in Meta-lemma 1, that we cannot introduce a finite number of arbitrary primitive recursive functions and relations, as function letters and predicate letters respectively, into a formal system of Peano Arithmetic, such as PA, without risking inconsistency. Ipso facto, PA+PRA is inconsistent.

It follows that, if the appended arguments are valid, then:

*f*) ZF is inconsistent.

## 2. Appendix

**Meta-theorem 1:** There is a primitive recursive relation that is not the standard interpretation of any of its formal representations in a formal Peano Arithmetic PA.

**Proof:** We consider Gödel's primitive recursive relation $\sim xB_{PA}(Sb(y\ 19|Z(y)))$.

(*a*) We assume that every primitive recursive function or relation is the standard interpretation of at least one of its formal representations in PA.

(*b*) There is, thus, some PA-formula $[\sim xB_{PA}(Sb(y\ 19|Z(y)))]$ whose standard interpretation is the primitive recursive relation $\sim xB_{PA}(Sb(y\ 19|Z(y)))$.

(*c*) Now, in every model M (cf. [Me64], p192-3) of PA, we can also interpret:

(*i*) the integer 0 as the interpretation of the symbol "0";

(*ii*) the successor operation as the interpretation of the successor function "";

(*iii*) addition and multiplication as the interpretations of "+" and ".";

(*iv*) the interpretation of the predicate letter = as the identity relation.

(*d*) Since the numerals of PA interpret as a sub-domain of every model M of PA, the natural numbers are, then, a sub-domain of every M.

(*e*) Further, by the hypothesis (*a*), all of Gödel's 45 primitive recursive functions and relations ([Go31a], p17-22) are also, then, mirrored in every model M of PA, and the PA-formula, $[\sim xB_{PA}(Sb(p\ 19|Z(p)))]$, always interprets as the M-relation $\sim xB_{PA}(Sb(p\ 19|Z(p)))$, where $[p]$ is the numeral that represents the natural number $p$ in PA, and $p$ is the Gödel-number of the PA-formula $[(Ax)\sim xB_{PA}(Sb(y\ 19|Z(y)))]$.

(*f*) Hence, in every model M of PA, the relation $\sim xB_{PA}(Sb(p\ 19|Z(p)))$ holds in M if, and only if, $x$ is a M-number that is not the Gödel-number of a proof of $[(Ax)\sim xB_{PA}(Sb(p\ 19|Z(p)))]$ in PA.

(*g*) Further, since the Gödel-number of a proof of [(A*x*)~*xB*$_{PA}$(*Sb*(*p* 19|*Z*(*p*)))] in PA is necessarily a natural number, ~*xB*$_{PA}$(*Sb*(*p* 19|*Z*(*p*))) holds in every model M of PA if *x* is an M-number that is not a natural number.

(*h*) Now, Gödel has shown (cf. [Go31a], Theorem VI) that ~*xB*$_{PA}$(*Sb*(*p* 19|*Z*(*p*))) holds over the domain of the natural numbers.

(*i*) It follows that ~*xB*$_{PA}$(*Sb*(*p* 19|*Z*(*p*))) is satisfied by all *x* in every model M of PA.

(*j*) Hence, the PA-formula [(A*x*)~*xB*$_{PA}$(*Sb*(*p* 19|*Z*(*p*)))] is true (cf. [Me64], p51) in every model M of PA, and, by Gödel's Completeness Theorem ([Me64], Corollary 2.14, p68), [(A*x*)~*xB*$_{PA}$(*Sb*(*p* 19|*Z*(*p*)))] is PA-provable.

However, since Gödel has shown that [(A*x*)~*xB*$_{PA}$(*Sb*(*p* 19|*Z*(*p*)))] is not PA-provable (cf. [Go31a], Theorem VI), we conclude that assumption (*a*) does not hold. This proves the theorem.¶

**Meta-lemma 1:** We cannot introduce a finite number of arbitrary recursive number-theoretic functions and relations, as function letters and predicate letters respectively, into PA without risking inconsistency.

**Proof:** Adding "*B*$_{PA}$", and a finite number of other functions and relations in terms of which it is defined (cf. [Go31a], p17-22, def. 1-45), as new function letters and predicate letters, respectively, to PA, along with associated defining axioms (cf. [Me64], §9, p82), would yield a formal system PA* in which, by the arguments of Meta-theorem 1, the PA*-formulas, [(A*x*)~*xB*$_{PA}$(*Sb*(*p* 19|*Z*(*p*)))], and therefore [~*xB*$_{PA}$(*Sb*(*p* 19|*Z*(*p*)))][2], are true in every model of PA*, and hence PA*-provable.

---

[2] By Generalisation

(Note that the argument is *not* that the PA* formula, $[(Ax) \sim xB_{PA*}(Sb(p\ 19|Z(p)))]$, is PA*-provable.)

Since all the added functions and relations are primitive recursive, they are strongly represented in PA ([Me64], Proposition 3.23, p134). Every proof sequence of PA* can, thus, be converted into a proof sequence of PA ([Me64], Proposition 2.29, pp82-83). Ergo, if Gödel's primitive recursive relation, $\sim xB_{PA}(Sb(p\ 19|Z(p)))$, is represented in PA by $[G(x)]$, then both $[G(x)]$, and hence $[(Ax)G(x)]$, are PA-provable - contradicting Gödel's Theorem VI (cf. [Go31a], p24). This proves the lemma.¶

*(Created: Friday 18th February 2005 9:41:53 AM IST by re@alixcomsi.com. Last updated: Sunday 27th February 2005 12:07:57 AM IST by re@alixcomsi.com)*